\input amstex
\magnification=\magstep0
\documentstyle{amsppt}
\pagewidth{6.0in}
\input amstex
\topmatter
\title
Explicit computations with the Divided Symmetrization operator
\endtitle
\author
Tewodros Amdeberhan\\
Tulane University \\
tamdeber\@tulane.edu
\endauthor
\abstract
Given a multi-variable polynomial, there is an associated divided symmetrization (in particular turning it into a symmetric function). Postinkov has found the volume of a permutohedron as a divided symmetrization (DS) of the power of a certain linear form. The main task in this paper is to exhibit and prove closed form DS-formulas for a variety of polynomials. We hope the results to be valuable and available to the research practitioner in these areas. Also, the methods of proof utilized here are simple and amenable to many more analogous computations. We conclude the paper with a list of such formulas.
\endabstract
\endtopmatter
\def\({\left(}
\def\){\right)}

\document

\noindent
Throughout, let $S_n$ denote the symmetric group of permutations of the $n$-element set $\{1,2,\dots,n\}$. Given a function $f(\lambda_1,\dots,\lambda_n)$, the \it divided symmetrization \rm (DS) of $f$, is defined by
$$\left\langle f\right\rangle=\sum_{\sigma\in S_n}f(\lambda_{\sigma(1)},\dots,\lambda_{\sigma(n)})
\prod_{k=1}^{n-1}\frac1{\lambda_{\sigma(k)}-\lambda_{\sigma(k+1)}}.$$
Our motivation comes from a beautiful work [1] of Alexander Postnikov who has found the volume of a permutohedron in terms of divided symmetrization of a certain expression, in addition he offers a combinatorial interpretation of the resulting coefficients. To put in context, we recall the following:
\smallskip
\noindent
\bf Proposition 0 [Postnikov]. \it If $f$ is a polynomial of degree $n-1$ in the variables $\lambda_1,\dots,\lambda_n$, then
its divided symmetrization $\langle f\rangle$ is a constant. If $\deg f < n-1$, then $\langle f\rangle=0$. \rm
\smallskip
\noindent
\bf Proof. \rm Write $\langle f\rangle=g/\Delta$ where $\Delta=\prod_{i<j}(\lambda_i-\lambda_j)$ is the antisymmetric Vandermonde
determinant and $g$ is a some antisymmetric polynomial. Since $\langle f\rangle$ is symmetric, $g$ is divisible by $\Delta$. Because
$\deg g=\deg\Delta=\binom{n}2$, their quotient must be a constant. A similar argument shows  if $\deg g<n-1$ then $\deg g<\deg\Delta$
and, hence, $g=0$. The proof is complete. $\square$
\bigskip
\noindent
This proposition is applicable to (most) identities discussed in this note. But then, we ask: 

\it  is it possible to obtain the actual  values of these constants? \rm

\noindent
If so, then: 

\it do these constants have a closed form? \rm

\noindent
In this paper, we answer both questions for a variety of polynomial functions. The selection of the results is not meant to be representative, rather reflects the author's personal taste. For starters, we invite the reader to ponder on proving the following three numerical identities before continuing further on.

\pagebreak

\smallskip
\noindent
\it Involving products\rm:
$$\sum_{\sigma\in S_{n+1}}\prod_{k=1}^n\frac{\sigma(1)+\cdots+\sigma(k)}{\sigma(k)-\sigma(k+1)}=n!.$$
$$\sum_{\sigma\in S_{n+1}}\sigma(1)^m\prod_{k=1}^n\frac{\sigma(1)-\sigma(k+1)}{\sigma(k)-\sigma(k+1)}
=1^m+2^m+\cdots+(n+1)^m.$$

\noindent
\it Involving sums\rm:
$$\align
\sum_{\sigma\in S_{n+1}}\sum_{k=1}^n
\frac{\sigma(1)+\cdots+\sigma(k)}{\sigma(k)-\sigma(k+1)}
&=\binom{n+1}2n!,\endalign$$
which equals the numbers of edges in the Hasse diagram of weak Bruhat order of $S_{n+1}$.

\bigskip
\noindent
From here on, as promised, we state and prove numerous divided difference formulas.
\bigskip
\noindent
\bf Lemma 1. \it If $\lambda_1,\lambda_2,\dots$ are variables, then
$$\sum_{\sigma\in S_{n+1}}\sum_{k=1}^n\frac{\lambda_{\sigma(k)}}{\lambda_{\sigma(k)}-\lambda_{\sigma(k+1)}}
=\binom{n+1}2n!.$$ \rm
\bf Proof. \rm Break up the sum according to $\sigma(k)=i, \sigma(k+1)=j$ so that
$$\align
\sum_{\sigma\in S_{n+1}}\sum_{k=1}^n\frac{\lambda_{\sigma(k)}}{\lambda_{\sigma(k)}-\lambda_{\sigma(k+1)}}
&=\sum_{k=1}^n\sum\Sb i,j=1 \\ i\neq j\endSb^{n+1} \frac{\lambda_i}{\lambda_i-\lambda_j}\sum_{\sigma'\in S_{n-1}}1
=(n-1)!\sum_{k=1}^n\sum\Sb i,j=1 \\ i< j\endSb^{n+1} \frac{\lambda_i-\lambda_j}{\lambda_i-\lambda_j} \\
&=(n-1)!\cdot\binom{n+1}2\sum_{k=1}^n1=\binom{n+1}2n!. \qquad \qquad \qquad \square
\endalign$$
\bf Corollary 2. \it If $\lambda_1,\lambda_2,\dots$ are variables, then
$$\sum_{\sigma\in S_{n+1}}\sum_{k=1}^n\frac{\lambda_{\sigma(1)}+\cdots+\lambda_{\sigma(k)}}{\lambda_{\sigma(k)}-\lambda_{\sigma(k+1)}}=\binom{n+1}2n!.$$ \rm
\bf Proof. \rm If  $j\neq k, k+1$ then
$\sum_{\sigma\in S_{n+1}}\frac{\lambda_{\sigma(j)}}{\lambda_{\sigma(k)}-\lambda_{\sigma(k+1)}}=0$. Thus by Lemma 1,
$$\sum_{\sigma\in S_{n+1}}\sum_{k=1}^n\frac{\lambda_{\sigma(1)}+\cdots+\lambda_{\sigma(k)}}{\lambda_{\sigma(k)}-\lambda_{\sigma(k+1)}}=
\sum_{\sigma\in S_{n+1}}\sum_{k=1}^n\frac{\lambda_{\sigma(k)}}{\lambda_{\sigma(k)}-\lambda_{\sigma(k+1)}}
=\binom{n+1}2n!. \qquad \qquad \square$$
\bf Lemma 3. \it If $\lambda_1,\lambda_2,\dots$ are variables, then
$$\sum_{\sigma\in S_{n+1}}\sum_{k=1}^n\frac{\lambda_{\sigma(1)}^2+\cdots+\lambda_{\sigma(k)}^2}{\lambda_{\sigma(k)}-\lambda_{\sigma(k+1)}}=n(\lambda_1+\cdots+\lambda_{n+1})n!.$$ \rm
\bf Proof. \rm Since the terms $\lambda_{\sigma(\ell)}^2$ contribute to zero, whenever $\ell<k$, we restrict to $\ell=k$ and proceed as in the proof of Lemma 1 with $\sigma(k)=i, \sigma(k+1)=j$. That means,
$$\align
\sum_{\sigma\in S_{n+1}}\sum_{k=1}^n\frac{\lambda_{\sigma(k)}^2}{\lambda_{\sigma(k)}-\lambda_{\sigma(k+1)}}
&=\sum_{k=1}^n\sum\Sb i,j=1 \\ i< j\endSb^{n+1} \frac{\lambda_i^2-\lambda_j^2}{\lambda_i-\lambda_j}\sum_{\sigma'\in S_{n-1}}1
=(n-1)!\sum\Sb i,j=1 \\ i< j\endSb^{n+1} (\lambda_i+\lambda_j)\sum_{k=1}^n1\\
&=n!\sum\Sb i,j=1 \\ i< j\endSb^{n+1} (\lambda_i+\lambda_j)
=n!\cdot n(\lambda_1+\cdots+\lambda_{n+1}).  \qquad \qquad \square \endalign$$

\pagebreak

\noindent
\bf Lemma 4. \it If $Y, \lambda_1,\lambda_2,\dots$ be variables, then
$$\sum_{\sigma\in S_{n+1}}\prod_{k=1}^n
\frac{\lambda_{\sigma(k)}-Y}{\lambda_{\sigma(k)}-\lambda_{\sigma(k+1)}}=1 \qquad\text{and} \qquad
\sum_{\sigma\in S_{n+1}}\prod_{k=1}^n
\frac{Y-\lambda_{\sigma(k+1)}}{\lambda_{\sigma(k)}-\lambda_{\sigma(k+1)}}=1.$$ \rm

\noindent
\bf Proof. \rm We prove the first identity only, the same argument work for the second one. When $n=1$, the vacuous product confirms the assertion.  We induct on $n$. The base case $n=1$ is easily verifiable too, since $\frac{\lambda_1-Y}{\lambda_1-\lambda_2}+\frac{\lambda_2-Y}{\lambda_2-\lambda_1}=1$. Suppose the statement holds for $n-1$. For next step $n$, observe that the polynomial on the left-hand side is of degree at most $n$, in $Y$. Let's compute its values at the $n+1$ points $Y\in\{\lambda_1,\dots,\lambda_{n+1}\}$ as follows: for
each $1\leq j\leq n+1$,
$$\align \sum_{\sigma\in S_{n+1}}\prod_{k=1}^n
\frac{\lambda_{\sigma(k)}-\lambda_j}{\lambda_{\sigma(k)}-\lambda_{\sigma(k+1)}}&=
\sum_{\sigma'\in S_n}\frac{\lambda_{\sigma'(n)}-\lambda_j}{\lambda_{\sigma'(n)}-\lambda_{\sigma(n+1)}}
\prod_{k=1}^{n-1}\frac{\lambda_{\sigma'(k)}-\lambda_j}{\lambda_{\sigma'(k)}-\lambda_{\sigma'(k+1)}} \\
&=\sum_{\sigma'\in S_n}\frac{\lambda_{\sigma'(n)}-\lambda_j}{\lambda_{\sigma'(n)}-\lambda_j}
\prod_{k=1}^{n-1}\frac{\lambda_{\sigma'(k)}-\lambda_j}{\lambda_{\sigma'(k)}-\lambda_{\sigma'(k+1)}} \\
&=\sum_{\sigma'\in S_n}\prod_{k=1}^{n-1}\frac{\lambda_{\sigma'(k)}-\lambda_j}{\lambda_{\sigma'(k)}-\lambda_{\sigma'(k+1)}}=1;\endalign$$
where we used the induction assumption. So, the polynomial is a constant and the proof follows. $\square$
\bigskip
\noindent
\bf Lemma 5. \it As a direct consequence, we obtain the following partial fraction decomposition
$$\sum_{\sigma\in S_n}\frac1{\lambda_{\sigma(j)}-Y}\prod_{k=1}^{n-1}
\frac1{\lambda_{\sigma(k)}-\lambda_{\sigma(k+1)}}
=(-1)^{n-j}\binom{n-1}{j-1}\prod_{k=1}^n\frac1{\lambda_k-Y}; \qquad j\in\{1,\dots,n\}.$$ \rm
\bf proof. \rm We rewrite the assertion into an equivalent formulation and show that
$$\sum_{\sigma\in S_{n+1}}\frac{\lambda_{\sigma(n+1)}-Y}{\lambda_{\sigma(j)}-Y}\prod_{k=1}^n
\frac{\lambda_{\sigma(k)}-Y}{\lambda_{\sigma(k)}-\lambda_{\sigma(k+1)}}
=(-1)^{n-j+1}\binom{n}{j-1}.$$
Denote the sum on the left-hand side by $\Psi_{n,j}(Y)$. The case $j=n+1$ is simply the content of Lemma 4. So, assume $\sigma(n+1)\neq\sigma(j)$ (or $j\leq n$). By dropping-off vanishing terms along the way, we compute at $Y=\lambda_{\ell}$:
$$\align \Psi_{n,j}(\lambda_{\ell})
&=\sum_{S_{n+1}}\prod_{k=j}^n\frac{\lambda_{\sigma(k+1)}-\lambda_{\ell}}{\lambda_{\sigma(k)}-\lambda_{\sigma(k+1)}}\prod_{k=1}^{j-1}\frac{\lambda_{\sigma(k)}-\lambda_{\ell}}{\lambda_{\sigma(k)}-\lambda_{\sigma(k+1)}} \\
&=-\sum\Sb\sigma\in S_{n+1}\\ \sigma(j)=\ell\endSb
\prod_{k=j+1}^n\frac{\lambda_{\sigma(k+1)}-\lambda_{\ell}}{\lambda_{\sigma(k)}-\lambda_{\sigma(k+1)}}\prod_{k=1}^{j-2}\frac{\lambda_{\sigma(k)}-\lambda_{\ell}}{\lambda_{\sigma(k)}-\lambda_{\sigma(k+1)}}
\endalign$$
Now, form a disjoint partition of the underlying set of permutations as
$$\{\sigma\in S_{n+1}:\sigma(j)=\ell\}=\biguplus\Sb A\subset\{1,\dots,n+1\}-\{\ell\} \\
\text{$(j-1)$-subset}\endSb
\left\{\sigma\in S_{n+1}\vert \text{$\sigma(j)=\ell$ and $\{\sigma(1),\dots,\sigma(j-1)\}=A$}\right\}.$$
With this set-up and by Lemma 4, the above summation becomes
$$\align \Psi_{n,j}(\lambda_{\ell})&=-\sum\Sb A\subset\{1,\dots,n+1\}-\{\ell\} \\ \text{$(j-1)$-subset}\endSb
\sum_{S_{n-j+1}}\prod_{k=1}^{n-j}\frac{\lambda_{\sigma(k+1)}-\lambda_{\ell}}{\lambda_{\sigma(k)}-\lambda_{\sigma(k+1)}}\sum_{S_{j-1}}
\prod_{k=1}^{j-2}\frac{\lambda_{\sigma(k)}-\lambda_{\ell}}{\lambda_{\sigma(k)}-\lambda_{\sigma(k+1)}} \\
&=-\sum\Sb A\subset\{1,\dots,n+1\}-\{\ell\} \\ \text{$(j-1)$-subset}\endSb(-1)^{n-j}\cdot1
=(-1)^{n-j+1}\binom{n}{j-1}.
\endalign$$
The polynomial $\Psi_{n,j}(Y)$ being of degree at most $n$ and attaining the same value at the $n+1$ points $\lambda_1,\dots,\lambda_{n+1}$, it must be constant. The proof follows. $\square$

\bigskip
\noindent
\bf Lemma 6. \it If $Y, \lambda_1, \lambda_2, \dots$, then
$$\sum_{j=1}^n\prod\Sb k=1 \\ k\neq j\endSb^n \frac{Y-\lambda_j}{\lambda_k-\lambda_j}=1.$$ \rm
\bf Proof. \rm Induction on $n$, computing at $Y=\lambda_1, \dots, \lambda_n$. $\square$
\bigskip
\noindent
\bf Lemma 7. \it If $Y, \lambda_1, \lambda_2, \dots$ and $1\leq j\leq n+1$, then
$$\sum_{\sigma\in S_{n+1}}\prod_{k=1}^n\frac{\lambda_{\sigma(j)}-Y}{\lambda_{\sigma(k)}-\lambda_{\sigma(k+1)}}
=(-1)^{j-1}\binom{n}{j-1}.$$ \rm
\bf Proof. \rm Let $j=1$. From Lemma 5 and Lemma 6 (in that order), we obtain
$$\align
\sum_{\sigma\in S_{n+1}}\prod_{k=1}^n\frac{\lambda_{\sigma(j)}-Y}{\lambda_{\sigma(k)}-\lambda_{\sigma(k+1)}}
&=\sum_{\ell=1}^{n+1}
\sum_{\sigma'\in S_n}\frac{(\lambda_{\ell}-Y)^n}{\lambda_{\ell}-\lambda_{\sigma'(2)}}
\prod_{k=2}^n\frac1{\lambda_{\sigma'(k)}-\lambda_{\sigma'(k+1)}} \\
&=-(-1)^{n-1}\sum_{\ell=1}^{n+1}\prod\Sb k=1\\ k\neq\ell\endSb^{n+1}\frac{\lambda_{\ell}-Y}{\lambda_k-\lambda_{\ell}}=(-1)^n(-1)^n=1.\endalign$$
For the case $j=n+1$, once more, Lemma 5 and Lemma 6 yield
$$\align
\sum_{\sigma\in S_{n+1}}\prod_{k=1}^n\frac{\lambda_{\sigma(j)}-Y}{\lambda_{\sigma(k)}-\lambda_{\sigma(k+1)}}
&=\sum_{\ell=1}^{n+1}
\sum_{\sigma'\in S_n}\frac{(\lambda_{\ell}-Y)^n}{\lambda_{\sigma'(n)}-\lambda_{\ell}}
\prod_{k=1}^{n-1}\frac1{\lambda_{\sigma'(k)}-\lambda_{\sigma'(k+1)}} \\
&=\sum_{\ell=1}^{n+1}\prod\Sb k=1\\ k\neq\ell\endSb^{n+1}\frac{\lambda_{\ell}-Y}{\lambda_k-\lambda_{\ell}}
=(-1)^n.\endalign$$

\noindent
Suppose that $2\leq j\leq n$. By Lemma 5, Lemma 6 and Pascal's recurrence we find that
$$\align
\sum_{\sigma\in S_{n+1}}
\prod_{k=1}^n\frac{\lambda_{\sigma(j)}-Y}{\lambda_{\sigma(k)}-\lambda_{\sigma(k+1)}}
&=\sum\Sb\ell=1\endSb^{n+1}\sum_{\sigma'\in S_n}\frac{(\lambda_{\ell}-Y)^n}{(\lambda_{\sigma'(j-1)}-\lambda_{\ell})(\lambda_{\ell}-\lambda_{\sigma'(j+1)})}
\prod\Sb k=1\\ k\neq j, j-1\endSb^n\frac1{\lambda_{\sigma'(k)}-\lambda_{\sigma'(k+1)}} \\
&=\sum\Sb\ell=1\endSb^{n+1}\sum_{\sigma'\in S_n}\left[
\frac{(\lambda_{\ell}-Y)^n}{\lambda_{\sigma'(j-1)}-\lambda_{\ell}}\right]
\frac1{\lambda_{\sigma'(j-1)}-\lambda_{\sigma'(j+1)}}
\prod\Sb k=1\\ k\neq j-1, j\endSb^n\frac1{\lambda_{\sigma'(k)}-\lambda_{\sigma'(k+1)}} \\
&+\sum\Sb\ell=1\endSb^{n+1}\sum_{\sigma'\in S_n}\left[\frac{(\lambda_{\ell}-Y)^n}{\lambda_{\ell}-\lambda_{\sigma'(j+1)}}\right]
\frac1{\lambda_{\sigma'(j-1)}-\lambda_{\sigma'(j+1)}}
\prod\Sb k=1\\ k\neq j-1,j\endSb^n\frac1{\lambda_{\sigma'(k)}-\lambda_{\sigma'(k+1)}} \\
&=\left[(-1)^{n-j+1}\binom{n-1}{j-2}-(-1)^{n-j}\binom{n-1}{j-1}\right]
\sum_{\ell=1}^{n+1}
\prod\Sb k=1\\ k\neq \ell\endSb^{n+1}\frac{\lambda_{\ell}-Y}{\lambda_k-\lambda_{\ell}} \\
&=\left[(-1)^{n-j+1}\binom{n}{j-1}\right](-1)^n=(-1)^{j-1}\binom{n}{j-1}. \endalign$$
The proof is complete. $\square$

\bigskip
\noindent
\bf Corollary 8. \it If $\lambda_1,\lambda_2,\dots$ are variables and $e_n(\lambda)$ the $n^{th}$-elementary symmetric function, then
$$\sum_{\sigma\in S_{n+1}}\prod_{k=1}^n\frac{\lambda_{\sigma(k)}^2}{\lambda_{\sigma(k)}-\lambda_{\sigma(k+1)}}
=e_n(\lambda_1,\dots,\lambda_{n+1}).$$ \rm
\bf Proof. \rm We prove the equivalent claim
$\sum_{\sigma\in S_{n+1}}\frac1{\lambda_{\sigma(n+1)}}\prod_{k=1}^n\frac{\lambda_{\sigma(k)}}{\lambda_{\sigma(k)}-\lambda_{\sigma(k+1)}}=\sum_{j=1}^{n+1}\frac1{\lambda_j}$. To this end,
$$\align \sum_{S_{n+1}}\frac1{\lambda_{\sigma(n+1)}}
\prod_{k=1}^n\frac{\lambda_{\sigma(k)}}{\lambda_{\sigma(k)}-\lambda_{\sigma(k+1)}}
&=\sum_{S_{n+1}}\frac{\lambda_{\sigma(n)}}{\lambda_{\sigma(n+1)}(\lambda_{\sigma(n)}-\lambda_{\sigma(n+1)})}
\prod_{k=1}^{n-1}\frac{\lambda_{\sigma(k)}}{\lambda_{\sigma(k)}-\lambda_{\sigma(k+1)}} \\
&=\sum_{S_{n+1}}\left[\frac1{\lambda_{\sigma(n+1)}}+\frac1{\lambda_{\sigma(n)}-\lambda_{\sigma(n+1)}}\right]
\prod_{k=1}^{n-1}\frac{\lambda_{\sigma(k)}}{\lambda_{\sigma(k)}-\lambda_{\sigma(k+1)}}.
\endalign$$
Next, we separate the last sum into two and apply Lemma 4 to the first summation as follows:
$$\align \sum_{S_{n+1}}\frac1{\lambda_{\sigma(n+1)}}
\prod_{k=1}^{n-1}\frac{\lambda_{\sigma(k)}}{\lambda_{\sigma(k)}-\lambda_{\sigma(k+1)}}
&=\sum_{j=1}^{n+1}\frac1{\lambda_j}\sum_{S_n}\prod_{k=1}^{n-1}\frac{\lambda_{\sigma(k)}}{\lambda_{\sigma(k)}-\lambda_{\sigma(k+1)}}=\sum_{j=1}^{n+1}\frac1{\lambda_j}.\endalign$$
Since $f=\prod_{k=1}^{n-1}\lambda_k$ is of degree $n-1$, Proposition 0 implies the vanishing of the second sum:
$$\align \sum_{S_{n+1}}\frac1{\lambda_{\sigma(n)}-\lambda_{\sigma(n+1)}}
\prod_{k=1}^{n-1}\frac{\lambda_{\sigma(k)}}{\lambda_{\sigma(k)}-\lambda_{\sigma(k+1)}}
&=\sum_{S_{n-1}}\prod_{k=1}^{n-1}\lambda_{\sigma(k)}
\prod_{k=1}^n\frac1{\lambda_{\sigma(k)}-\lambda_{\sigma(k+1)}}=0.\endalign$$
The proof is complete. $\square$

\bigskip
\noindent
\bf Corollary 9. \it If $Y, \lambda_1,\lambda_2,\dots$ are variables, then
$$\sum_{\sigma\in S_{n+1}}\prod_{k=1}^n
\frac{\lambda_{\sigma(1)}-\lambda_{\sigma(k+1)}-Y}{\lambda_{\sigma(k)}-\lambda_{\sigma(k+1)}}=n+1.$$ \rm
\bf Proof. \rm Rewrite the sum on the left-hand side and apply Lemma 5 as follows:
$$\align
-\frac1{Y}\sum_{S_{n+1}}\frac{\prod_{k=1}^{n+1}\lambda_{\sigma(1)}-\lambda_k-Y}
{\prod_{k=1}^n\lambda_{\sigma(k)}-\lambda_{\sigma(k+1)}}
&=\frac1{Y}\sum_{j=1}^{n+1}\prod_{k=1}^{n+1}(\lambda_j-\lambda_k-Y)
\sum_{S_n}\frac1{\lambda_{\sigma(2)}-\lambda_j}\prod_{k=2}^n\frac1{\lambda_{\sigma(k)}-\lambda_{\sigma(k+1)}} \\
&=\frac1{Y}\sum_{j=1}^{n+1}\prod_{k=1}^{n+1}(\lambda_j-\lambda_k-Y)\left[(-1)^{n-1}
\prod\Sb k=1\\ k\neq j\endSb^{n+1}\frac1{\lambda_k-\lambda_j}\right] \\
&=\sum_{j=1}^{n+1}\prod\Sb k=1\\ k\neq j\endSb^{n+1}\frac{\lambda_k-\lambda_j+Y}{\lambda_k-\lambda_j}. \endalign$$
Since Proposition 0 ensures that the above sum is independent of both $Y$ and the $\lambda$'s; a choice of $Y=0$ settles the claim. The proof is complete. $\square$
\smallskip
\noindent
\bf Corollary 10. \it For an indeterminate $Y$, we have $\sum_{j=0}^n(-1)^j\binom{Y-1}{j}\binom{Y+n-j}{n-j}=n+1$. \rm
\smallskip
\noindent
\bf Proof. \rm In the proof of Corollary 9, choose $\lambda_k=k$ which is admissible due to Proposition 0. Thus,
$$n+1=\sum_{j=1}^{n+1}\prod\Sb k=1\\ k\neq j\endSb^{n+1}\frac{k-j+Y}{k-j}
=\sum_{j=0}^n(-1)^j\binom{Y-1}{j}\binom{Y+n-j}{n-j}.$$
Alternatively, this identity can be generalized to $\sum_j(-1)^j\binom{Y-k}j\binom{Y+n-j}{n-j}=\binom{n+k}k$ and proved by convolving generating functions as $(1-x)^{Y-k}(1-x)^{-Y-1}=(1-x)^{-k-1}$. The proof follows. $\square$
\smallskip
\noindent
\bf Corollary 11. \it If $Y, \lambda_1,\lambda_2,\dots$ are variables, then
$$\sum_{\sigma\in S_{n+1}}\prod_{k=1}^n
\frac{\lambda_{\sigma(1)}+\lambda_{\sigma(2)}-\lambda_{\sigma(k+1)}-Y}{\lambda_{\sigma(k)}-\lambda_{\sigma(k+1)}}=2^{n+1}-n-2.$$ \rm
\bf Proof. \rm Rearrange the left-hand side (LHS, for short), apply Lemma 5 and proceed as:
$$\align
LHS &=\sum\Sb i,j=1\\ i\neq j\endSb^{n+1}\frac{\prod_{k=1}^{n+1}\lambda_i+\lambda_j-\lambda_k-Y}{(\lambda_j-Y)(\lambda_i-\lambda_j)}\sum_{S_{n-1}}\frac1{\lambda_j-\lambda_{\sigma(3)}}\prod_{k=3}^n
\frac1{\lambda_{\sigma(k)}-\lambda_{\sigma(k+1)}} \\
&=(-1)^{n+1}\sum\Sb i,j=1\\ i\neq j\endSb^{n+1}
\frac{\prod_{k=1}^{n+1}\lambda_i+\lambda_j-\lambda_k-Y}{(\lambda_j-Y)(\lambda_i-\lambda_j)}
\prod\Sb k=1\\ k\neq i,j\endSb^{n+1}\frac1{\lambda_k-\lambda_j} \\
&=-\sum\Sb i,j=1\\ i\neq j\endSb^{n+1}\prod\Sb k=1\\ k\neq i\endSb^{n+1}(\lambda_k-\lambda_j+Y-\lambda_i)
\prod\Sb k=1\\ k\neq j\endSb^{n+1}(\lambda_k-\lambda_j)^{-1} \\
&=-\sum\Sb i,j=1\endSb^{n+1}\prod\Sb k=1\\ k\neq i\endSb^{n+1}(\lambda_k-\lambda_j+Y-\lambda_i)
\prod\Sb k=1\\ k\neq j\endSb^{n+1}(\lambda_k-\lambda_j)^{-1} +
\sum_{j=1}^{n+1}\prod\Sb k=1\\ k\neq j\endSb^{n+1}\frac{\lambda_k-2\lambda_j+Y}{\lambda_k-\lambda_j}\\
&=:\text{I$_{n}(Y)$+II$_{n+1}(Y)$}.\endalign$$
To determine the second sum, we induct on $n$. Let's evaluate II$_{n+1}(\lambda_{\ell})$ for each $1\leq\ell\leq n+1$:
$$\text{II}_{n+1}(\lambda_{\ell})=1+\sum\Sb j=1\\ j\neq\ell\endSb^{n+1}2
\prod\Sb k=1\\ k\neq j,\ell\endSb^{n+1}\frac{\lambda_k-2\lambda_j+\lambda_{\ell}}{\lambda_k-\lambda_j}=
1+2\text{II}_n(\lambda_{\ell}).$$
After consulting initial conditions, this recurrence implies II$_{n+1}(\lambda_{\ell})=2^{n+1}-1$. Since $\deg$ II$_{n+1}\leq n$, as a polynomial  in $Y$, it must be a constant. Next, we employ Proposition 0 which permits us working with $\lambda_k=k$ (even ignoring $Y$) and thus obtain
$$\align \text{I}_{n+1}(Y)
%&=\sum\Sb i,j=1\endSb^{n+1}(-1)^j\frac{j}{Y-j}\binom{Y+n+1-i-j}{n+1}\binom{n+1}j \\
&=\sum\Sb i,j=1\endSb^{n+1}(-1)^{n+1-j}\binom{i+j-1}{n+1}\binom{n+1}j
=\sum_{j=1}^{n+1}\binom{n+1}j(-1)^{n+1-j}\sum_{i=1}^{n+1}\binom{i+j-1}{n+1} \\
&=\sum_{j=1}^{n+1}(-1)^{n+1-j}\binom{n+1}j\sum_{k=0}^{j-1}\binom{n+1+k}k
=\sum_{j=1}^{n+1}(-1)^{n+1-j}\binom{n+1}j\binom{n+j+1}{j-1}=n+1. \endalign$$
The last identity is a consequence of the expansions of $(1-x)^{-Z-3}(1-x)^{Z+1}=(1-x)^{-2}$ which yield $\sum_{j=0}^n(-1)^{n-j}\binom{Z+j+2}j\binom{Z+1}{n-j}=n+1$ then simply replace $Z=n$ here.
In conclusion, we have -I$_{n+1}(Y)+$II$_{n+1}(Y)=-n-1+2^{n+1}-1$ as required. The proof is complete. $\square$

\bigskip
\noindent
\bf Corollary 12. \it Let $\lambda_1,\lambda_2,\dots$ be variables and $m\geq0$ an integer. Then, we have
$$\sum_{\sigma\in S_{n+1}}\lambda_{\sigma(1)}^m\prod_{k=1}^n\frac{\lambda_{\sigma(1)}-\lambda_{\sigma(k+1)}}{\lambda_{\sigma(k)}-\lambda_{\sigma(k+1)}}=\lambda_1^m+\cdots+\lambda_{n+1}^m.$$\rm
\bf Proof. \rm A simple rewrite together with Lemma 4 produce
$$\sum_{\sigma\in S_{n+1}}\lambda_{\sigma(1)}^m\prod_{k=1}^n\frac{\lambda_{\sigma(1)}-\lambda_{\sigma(k+1)}}{\lambda_{\sigma(k)}-\lambda_{\sigma(k+1)}}=\sum_{j=1}^{n+1}\lambda_j^m\sum_{S_n}
\prod_{k=2}^n\frac{\lambda_j-\lambda_{\sigma(k+1)}}{\lambda_{\sigma(k)}-\lambda_{\sigma(k+1)}}=
\sum_{j=1}^{n+1}\lambda_j^m. \qquad \square$$
\bf Lemma 13. \it Let $Y, Z, \lambda_1,x_1,\lambda_2,x_2,\dots$ be variables. Then, we have
$$\sum_{\sigma\in S_{n+1}}\frac1{\lambda_{\sigma(1)}-Z}\prod_{k=1}^n
\frac{\lambda_{\sigma(1)}x_1+\cdots+\lambda_{\sigma(k)}x_k-Y}{\lambda_{\sigma(k)}-\lambda_{\sigma(k+1)}}
=\prod_{k=1}^n(Y-(x_1+\cdots+x_k)Z)\prod_{k=1}^{n+1}(\lambda_k-Z)^{-1}.$$ \rm
\bf Proof. \rm We proceed by induction. Assume the statement is valid for $n$. So, for $n+1$, we find that
$$\align \sum_{S_{n+1}}\frac1{\lambda_{\sigma(1)}-Z}\prod_{k=1}^n
\frac{\sum_{i=1}^k\lambda_{\sigma(i)}x_i-Y}{\lambda_{\sigma(k)}-\lambda_{\sigma(k+1)}}
&=\sum_{j=1}^{n+1}\frac{\lambda_jx_1-Y}{\lambda_j-Z}\sum_{S_n}\frac{-1}{\lambda_{\sigma(2)}-\lambda_j}
\prod_{k=2}^n\frac{\sum_{i=2}^k\lambda_{\sigma(i)}x_i-(Y-\lambda_jx_1)}
{\lambda_{\sigma(k)}-\lambda_{\sigma(k+1)}} \\
&=\sum_{j=1}^{n+1}\frac{Y-\lambda_jx_1}{\lambda_j-Z}
\prod_{k=1}^{n-1}(Y-\lambda_jx_1-\sum_{i=1}^kx_{i+1}\lambda_j)\prod\Sb k=1\\ k\neq j\endSb^{n+1}(\lambda_k-\lambda_j)^{-1} \\
&=\sum_{j=1}^{n+1}\frac1{\lambda_j-Z}
\prod_{k=1}^n(Y-\sum_{i=1}^kx_i\lambda_j)\prod\Sb k=1\\ k\neq j\endSb^{n+1}(\lambda_k-\lambda_j)^{-1} \\
&=\prod_{k=1}^n(Y-\sum_{i=1}^kx_iZ)\prod_{k=1}^{n+1}(\lambda_k-Z)^{-1}. \endalign$$
The last equality results from the standard partial fraction decomposition applied to the rational function
$\prod_{k=1}^n(Y-\sum_{i=1}^kx_iZ)\prod_{k=1}^{n+1}(\lambda_k-Z)^{-1}$, in the variable $Z$. The proof complete. $\square$

\pagebreak

\bigskip
\noindent
\bf Corollary 14. \it Let $Y, \lambda_1,x_,\lambda_2,x_2,\dots$ be variables. Then, we have
$$\sum_{\sigma\in S_{n+1}}\prod_{k=1}^n
\frac{\lambda_{\sigma(1)}x_1+\cdots+\lambda_{\sigma(k)}x_k-Y}{\lambda_{\sigma(k)}-\lambda_{\sigma(k+1)}}
=\prod_{k=1}^n(x_1+\cdots+x_k).$$ \rm
\bf Proof. \rm A reduction of the sum to $S_n$ and a use of Lemma 13 (with $Z\rightarrow\lambda_j, Y\rightarrow Y-\lambda_jx_1$) yields
$$\align \sum_{S_{n+1}}\prod_{k=1}^n
\frac{\sum_{i=1}^k\lambda_{\sigma(i)}x_i-Y}{\lambda_{\sigma(k)}-\lambda_{\sigma(k+1)}}
&=\sum_{j=1}^{n+1}(Y-\lambda_jx_1)\sum_{S_n}\frac1{\lambda_{\sigma(2)}-\lambda_j}
\prod_{k=2}^n \frac{\sum_{i=2}^k\lambda_{\sigma(i)}x_i-(Y-\lambda_jx_1)}
{\lambda_{\sigma(k)}-\lambda_{\sigma(k+1)}} \\
&=\sum_{j=1}^{n+1}(Y-\lambda_jx_1)
\prod_{k=1}^{n-1}(Y-\lambda_jx_1-\sum_{i=1}^kx_{i+1}\lambda_j)\prod\Sb k=1\\ k\neq j\endSb^{n+1}(\lambda_k-\lambda_j)^{-1} \\
&=\sum_{j=1}^{n+1}\prod_{k=1}^n(Y-\sum_{i=1}^kx_i\lambda_j)\prod\Sb k=1\\ k\neq j\endSb^{n+1}(\lambda_k-\lambda_j)^{-1}=\prod_{k=1}^n(x_1+\cdots+x_k).
\endalign$$
To justify the last equality, rearrange the partial fraction seen in the proof of Lemma 13 as
$$\sum_{j=1}^{n+1}\frac{\lambda_{n+1}-Z}{\lambda_j-Z}
\prod_{k=1}^n(Y-\sum_{i=1}^kx_i\lambda_j)\prod\Sb k=1\\ k\neq j\endSb^{n+1}(\lambda_k-\lambda_j)^{-1}
=\prod_{k=1}^n\frac{Y-\sum_{i=1}^kx_iZ}{\lambda_k-Z},$$
and take the limit $Z\rightarrow\infty$ on both sides. The proof follows. $\square$
\bigskip
\noindent
\bf Corollary 15. \it Let $Y, \lambda_1,\lambda_2,\dots$ be variables. Then, we have
$$\sum_{\sigma\in S_{n+1}}\lambda_{\sigma(1)}\prod_{k=1}^n
\frac{\lambda_{\sigma(1)}+\cdots+\lambda_{\sigma(k)}}{\lambda_{\sigma(k)}-\lambda_{\sigma(k+1)}}
=n!(\lambda_1+\cdots+\lambda_{n+1}).$$ \rm
\bf Proof. \rm As usual, begin by reduction to $S_n$ and apply Lemma 13 (with $Z=\lambda_j, Y=-\lambda_j$) so that
$$\align \sum_{S_{n+1}}\lambda_{\sigma(1)}\prod_{k=1}^n
\frac{\sum_{i=1}^k\lambda_{\sigma(i)}}{\lambda_{\sigma(k)}-\lambda_{\sigma(k+1)}}
&=-\sum_{j=1}^{n+1}\lambda_j^2\sum_{S_n}\frac1{\lambda_{\sigma(2)}-\lambda_j}\prod_{k=2}^{n+1}
\frac{\sum_{i=2}^k\lambda_{\sigma(i)}+\lambda_j}{\lambda_{\sigma(k)}-\lambda_{\sigma(k+1)}} \\
&=-\sum_{j=1}^{n+1}\lambda_j^2\prod_{k=1}^{n-1}(-\lambda_j-k\lambda_j)\prod\Sb k=1\\ k\neq j\endSb^{n+1}
(\lambda_k-\lambda_j)^{-1} \\
&=n!\sum_{j=1}^{n+1}\lambda_j^{n+1}\prod\Sb k=1\\ k\neq j\endSb^{n+1}(\lambda_j-\lambda_k)^{-1}
=n!\sum_{k=1}^{n+1}\lambda_k. \endalign$$
We validate the last equality using Lemma 6 (with $Y=0$) and induction on $n$ based on the relation
$$\align \sum_{j=1}^{n+1}\lambda_j^{n+1}\prod\Sb k=1\\ k\neq j\endSb^{n+1}(\lambda_j-\lambda_k)^{-1}
&=\sum_{j=1}^n\frac{\lambda_j^n(\lambda_j-\lambda_{n+1}+\lambda_{n+1})}{\lambda_j-\lambda_{n+1}}
\prod\Sb k=1\\ k\neq j\endSb^n\frac1{\lambda_j-\lambda_k} +
\lambda_{n+1}^{n+1}\prod_{k=1}^n\frac1{\lambda_{n+1}-\lambda_k} \\
&=\sum_{j=1}^n\lambda_j^n\prod\Sb k=1\\ k\neq j\endSb^n(\lambda_j-\lambda_k)^{-1}+
\lambda_{n+1}\sum_{j=1}^{n+1}\prod\Sb k=1\\ k\neq j\endSb^{n+1}\frac{\lambda_j}{\lambda_j-\lambda_k} \\
&=\sum_{k=1}^n\lambda_k+\lambda_{n+1}. \qquad \qquad \square \endalign$$

\bigskip
\noindent
\bf Corollary 16. \it We have the identity (including its $q$-analogue)
$$\sum_{i=0}^n(-1)^{n-j}\binom{n}i(i+a)^{n+1}=n!\binom{n+1}2+(n+1)!a.$$ \rm
\bf Proof. \rm Although this result is well-known we present a different approach from a generalized fact that we have seen in the proof of Corollary 15, namely $n!\sum_{j=1}^{n+1}\lambda_j^{n+1}\prod\Sb k=1\\ k\neq j\endSb^{n+1}(\lambda_j-\lambda_k)^{-1}=n!\sum_{k=1}^{n+1}\lambda_k$. With the choice $\lambda_j=j+a-1$ and after some re-indexing the proof follows. $\square$
\smallskip
\noindent
\bf Corollary 17. \it We have the identity (including its $q$-analogue)
$$\sum_{i=0}^n(-1)^{n-i}\binom{n}i(i+a)^n=n!.$$ \rm
\bf Proof. \rm This, too, is a familiar identity. Yet, it follows from the more general relation that has been established in the
proof of Corollary 14; that is, $\sum_{j=1}^{n+1}\prod_{k=1}^n(Y-k\lambda_j)\prod\Sb k=1\\ k\neq j\endSb^{n+1}(\lambda_k-\lambda_j)^{-1}=n!$. Replace $Y=0$ and $\lambda_j=j+a-1$ to arrive at the desired conclusion. $\square$

\smallskip
\noindent
\bf Lemma 18. \it Let $Y, Z, \lambda_1,\lambda_2,\dots$ be variables. Then, we have
$$\sum_{\sigma\in S_{n+1}}\frac1{\lambda_{\sigma(n+1)}-Z}\prod_{k=1}^n
\frac{Y-\lambda_{\sigma(n+1)}+\cdots+\lambda_{\sigma(k+1)}}{\lambda_{\sigma(k)}-\lambda_{\sigma(k+1)}}
=\prod_{k=1}^n(Y-kZ)\prod_{k=1}^{n+1}(\lambda_k-Z)^{-1}.$$ \rm
\bf Proof. \rm This might be considered a "dual" version with a similar proof to that of Lemma 13. $\square$
\smallskip
\noindent
\bf Lemma 19. \it Let $Y, \lambda_1,\lambda_2,\dots$ be variables and $H_n:=\sum_{k=1}^n\frac1k$ the harmonic numbers. Then,
$$\sum_{\sigma\in S_{n+1}}\prod_{k=1}^n
\frac{\lambda_{\sigma(1)}+\cdots+\lambda_{\sigma(k)}-Y}
{(\lambda_{\sigma(k)}-Y)(\lambda_{\sigma(k)}-\lambda_{\sigma(k+1)})}=\frac{n!(1-H_n)(\lambda_1+\cdots+\lambda_{n+1}-Y)}{\prod_{k=1}^{n+1}\lambda_k-Y}.$$ \rm
\bf Proof. \rm We convert the left-hand side into an equivalent form and prove that
$$\sum_{S_{n+1}}(\lambda_{\sigma(n+1)}-Y)\prod_{k=1}^n
\frac{\sum_{i=1}^k\lambda_{\sigma(i)}-Y}{\lambda_{\sigma(k)}-\lambda_{\sigma(k+1)}}
=n!(1-H_n)\left(\sum_{k=1}^{n+1}\lambda_k-Y\right).$$
The second half of the last identity follows from Corollary 14, so we only consider the first half with the help of Lemma 18 (use $S_{n+1}\rightarrow S_n, Z\rightarrow\lambda_j, Y\rightarrow\lambda-Y-\lambda_j$) where $\vert\lambda\vert:=\sum_{i=1}^{n+1}\lambda_i$. Namely,
$$\align \sum_{S_{n+1}}\lambda_{\sigma(n+1)}\prod_{k=1}^n
\frac{\sum_{i=1}^k\lambda_{\sigma(i)}-Y}{\lambda_{\sigma(k)}-\lambda_{\sigma(k+1)}}
&=\sum_{S_{n+1}}\lambda_{\sigma(n+1)}\prod_{k=1}^n
\frac{(\vert\lambda\vert-Y)-\sum_{i=k+1}^{n+1}\lambda_{\sigma(i)}}{\lambda_{\sigma(k)}-\lambda_{\sigma(k+1)}} \\
&=\sum_{j=1}^{n+1}\lambda_j\sum_{S_n}\frac{(\vert\lambda\vert-Y)-\lambda_j}{\lambda_{\sigma(n)}-\lambda_j}\prod_{k=1}^{n-1}
\frac{(\vert\lambda\vert-Y-\lambda_j)-\sum_{i=k+1}^n\lambda_{\sigma(i)}}
{\lambda_{\sigma(k)}-\lambda_{\sigma(k+1)}} \\
&=\sum_{j=1}^{n+1}\lambda_j(\vert\lambda\vert-Y-\lambda_j)\prod_{k=1}^{n-1}(\vert\lambda\vert-Y-(k+1)\lambda_j)
\prod\Sb k=1\\ k\neq j\endSb^{n+1}(\lambda_k-\lambda_j)^{-1} \\
&=\sum_{j=1}^{n+1}\lambda_j\prod_{k=1}^n(\vert\lambda\vert-Y-k\lambda_j)
\prod\Sb k=1\\ k\neq j\endSb^{n+1}(\lambda_k-\lambda_j)^{-1}.\endalign$$
Expand the polynomial $P_n(Z):=Z\prod_{k=1}^n(\vert\lambda\vert-Y-kZ)+n!\prod_{k=1}^{n+1}(\lambda_k-Z)$, of degree $n$, in the Lagrange interpolating basis: $P_n(Z)=\sum_{j=1}^{n+1}\lambda_j
\prod_{k=1}^n(\vert\lambda\vert-Y-k\lambda_j)\prod\Sb k=1\\ k\neq j\endSb^{n+1}\frac{\lambda_k-Z}{\lambda_k-\lambda_j}$. Extracting the coefficients of $Z^n$ from both sides leads to
$n!(1-H_n)\vert\lambda\vert+n!H_nY$. The proof is complete. $\square$

\smallskip
\noindent
\bf Lemma 20. \it  Denote the Eulerian polynomials by $A_n(t)$. Let $\lambda_1,\lambda_2,\dots$ be variables.
Then, we have
$$\sum_{\sigma\in S_{n+1}}\prod_{k=1}^n\frac{\lambda_{\sigma(k)}-t\lambda_{\sigma(k+1)}}{\lambda_{\sigma(k)}-\lambda_{\sigma(k+1)}}=A_{n+1}(t)=\sum_{j=0}^nA(n+1,j)t^j.$$
\bf Proof. \rm Rewrite the left-hand side (denoted $F_{n+1}$) and apply Lemma 4 for a progressive expansion:
$$\align F_{n+1}
&=\sum_{S_{n+1}}\prod_{k=1}^n\left(1+\frac{(1-t)\lambda_{\sigma(k+1)}}{\lambda_{\sigma(k)}-\lambda_{\sigma(k+1)}}\right) \\
&=\sum_{S_{n+1}}\prod_{k=2}^n\left(1+\frac{(1-t)\lambda_{\sigma(k+1)}}{\lambda_{\sigma(k)}-\lambda_{\sigma(k+1)}}\right) +\sum_{S_{n+1}}\frac{(1-t)\lambda_{\sigma(2)}}{\lambda_{\sigma(1)}-\lambda_{\sigma(2)}}\prod_{k=2}^n
\left(1+\frac{(1-t)\lambda_{\sigma(k+1)}}{\lambda_{\sigma(k)}-\lambda_{\sigma(k+1)}}\right) \\
&=(n+1)F_n+\sum_{S_{n+1}}\frac{(1-t)\lambda_{\sigma(2)}}{\lambda_{\sigma(1)}-\lambda_{\sigma(2)}}\left[1+
\frac{(1-t)\lambda_{\sigma(3)}}{\lambda_{\sigma(2)}-\lambda_{\sigma(3)}}\right]\prod_{k=3}^n
\left(1+\frac{(1-t)\lambda_{\sigma(k+1)}}{\lambda_{\sigma(k)}-\lambda_{\sigma(k+1)}}\right) \\
&=(n+1)F_n+\binom{n+1}2\sum_{S_2}\frac{(1-t)\lambda_{\sigma(2)}}{\lambda_{\sigma(1)}-\lambda_{\sigma(2)}}\sum_{S_{n-1}}\prod_{k=3}^n
\left(1+\frac{(1-t)\lambda_{\sigma(k+1)}}{\lambda_{\sigma(k)}-\lambda_{\sigma(k+1)}}\right) \\
 & \qquad \qquad \qquad \qquad+\sum_{S_{n+1}}\frac{(1-t)^2\lambda_{\sigma(2)}\lambda_{\sigma(3)}}{(\lambda_{\sigma(1)}-\lambda_{\sigma(2)})(\lambda_{\sigma(2)}-\lambda_{\sigma(3)})}\prod_{k=3}^n
\left(1+\frac{(1-t)\lambda_{\sigma(k+1)}}{\lambda_{\sigma(k)}-\lambda_{\sigma(k+1)}}\right) \\
&=(n+1)F_n+\binom{n+1}2(t-1)F_{n-1}+\sum_{S_{n+1}}\prod_{k=1}^2\frac{(1-t)\lambda_{\sigma(k)}}{\lambda_{\sigma(k)}-\lambda_{\sigma(k+1)}}\prod_{k=4}^n
\left(1+\frac{(1-t)\lambda_{\sigma(k+1)}}{\lambda_{\sigma(k)}-\lambda_{\sigma(k+1)}}\right) \\
& \qquad \qquad \qquad \qquad+\sum_{S_{n+1}}\prod_{k=1}^3\frac{(1-t)\lambda_{\sigma(k)}}{\lambda_{\sigma(k)}-\lambda_{\sigma(k+1)}}\prod_{k=4}^n
\left(1+\frac{(1-t)\lambda_{\sigma(k+1)}}{\lambda_{\sigma(k)}-\lambda_{\sigma(k+1)}}\right) \\
&=\binom{n+1}1F_n+(t-1)\binom{n+1}2F_{n-1} \\
& \qquad \qquad \qquad \qquad +\binom{n+1}3\sum_{S_3}\prod_{k=1}^2\frac{(1-t)\lambda_{\sigma(k)}}{\lambda_{\sigma(k)}-\lambda_{\sigma(k+1)}}\sum_{S_{n-2}}\prod_{k=4}^n
\left(1+\frac{(1-t)\lambda_{\sigma(k+1)}}{\lambda_{\sigma(k)}-\lambda_{\sigma(k+1)}}\right) \\
& \qquad \qquad \qquad \qquad+\sum_{S_{n+1}}\prod_{k=1}^3\frac{(1-t)\lambda_{\sigma(k)}}{\lambda_{\sigma(k)}-\lambda_{\sigma(k+1)}}\prod_{k=4}^n
\left(1+\frac{(1-t)\lambda_{\sigma(k+1)}}{\lambda_{\sigma(k)}-\lambda_{\sigma(k+1)}}\right) \\
&=\binom{n+1}1F_n+(t-1)\binom{n+1}2F_{n-1}+(t-1)^2\binom{n+1}3F_{n-2}
+\sum_{S_{n+1}}\left(\prod_{k=1}^3\right)\left(\prod_{k=4}^n\right). \endalign$$
This process terminates with $F_{n+1}(t)=\sum_{j=0}^n(t-1)^{n-k}\binom{n+1}kF_{n+1-k}(t)$. So, the polynomials $F_n(t)$ fulfill the same recurrence as the
Eulerian polynomials $A_n(t)$ and also the same (check!) initial conditions. Hence $F_n=A_n$ and the proof is complete. $\square$

\bigskip
\noindent
\bf Corollary 21. \it Denote $[n]_q!:=\prod_{i=1}^n\frac{1-q^i}{1-q}$ and
$\left[n\atop j\right]:=\frac{[n]_q!}{[j]_q![n-j]_q!}$. We have the following identities
$$\align
\sum_{j=0}^n(-1)^j\binom{n}j(j+1)\prod_{k=1}^n\left[\binom{n+2}2-(j+1)k\right]&=n!^2(1-H_n)\binom{n+2}2, \\
\sum_{j=0}^n(-1)^j(1-q)q^{\binom{n+1}2+\binom{j+1}2+j}\left[ n \atop j \right]_q
%\left\{\aligned  n \\  j \endaligned\right \}_q
&=q^{\binom{n+1}2}(1-q)^{n+1}[n+1]_q!=\# GL(n+1,q);\endalign$$
where $GL(n,q)$ is the group of all $n\times n$ invertible matrices over the finite field $\Bbb{F}_q$ of $q$ elements. \rm
\smallskip
\noindent
\bf Proof. \rm At the end of the proof for Lemma 19, substitute $Y=0, \lambda_j=j$
(resp. $Y=\vert\lambda\vert, \lambda_j=q^j$) and simplify the sum to get the first (resp. the second) assertion. $\square$

\bigskip
\noindent
\bf Lemma 22. \it Let $Y, Z, \lambda_1,\lambda_2,\dots$ be variables. Then, we have
$$\sum_{\sigma\in S_{n+1}}\frac1{(\lambda_{\sigma(1)}-Y)(\lambda_{\sigma(n+1)}-Z)}\prod_{k=1}^n\frac1{\lambda_{\sigma(k)}-\lambda_{\sigma(k+1)}}=(Z-Y)^n\prod_{k=1}^{n+1}\frac1{(\lambda_k-Y)(\lambda_k-Z)}.$$ \rm
\bf Proof. \rm We proceed by induction on $n$. For $n+1$, the sum over $S_{n+1}$ is established by
$$\align \text{LHS}_{n+1}&=\sum\Sb i,j=1\\ i\neq j\endSb^{n+1}\frac{-1}{(\lambda_i-Y)(\lambda_j-Z)}\sum_{S_{n-1}}
\frac1{(\lambda_{\sigma(2)}-\lambda_i)(\lambda_{\sigma(n)}-\lambda_j)}\prod_{k=2}^{n-1}\frac1{\lambda_{\sigma(k)}-\lambda_{\sigma(k+1)}} \\
&=\sum\Sb i,j=1\\ i\neq j\endSb^{n+1}\frac{-1}{(\lambda_i-Y)(\lambda_j-Z)}(\lambda_j-\lambda_i)^{n-2}
\prod\Sb k=1\\ k\neq i,j\endSb^{n+1}\frac1{(\lambda_k-\lambda_i)(\lambda_k-\lambda_j)} \\
&=\sum_{i,j=1}^{n+1}\frac{(\lambda_j-\lambda_i)^n}{(\lambda_i-Y)(\lambda_j-Z)}
\prod\Sb k=1\\ k\neq i\endSb^{n+1}\frac1{\lambda_k-\lambda_i}
\prod\Sb k=1\\ k\neq j\endSb^{n+1}\frac1{\lambda_k-\lambda_j}
=(Z-Y)^n\prod_{k=1}^{n+1}\frac1{(\lambda_k-Y)(\lambda_k-Z)}.\endalign$$
The last equality is a consequence of partial fraction decompositions applied to the rational function
$(Z-Y)^n\prod_{k=1}^{n+1}
(\lambda_k-Y)^{-1}(\lambda_k-Z)^{-1}$, separately, in the variables $Y$ and $Z$. $\square$
%$$\sum_{\sigma\in S_{n+1}}\prod_{k=1}^{n-1}\frac{(\lambda_{\sigma(k+1)}-Y)(\lambda_{\sigma(k)}-Z)}
%{\lambda_{\sigma(k)}-\lambda_{\sigma(k+1)}}=(Z-Y)^n. \qquad \square$$
\bigskip
\noindent
\bf Lemma 23. \it Let $Y, \lambda_1,\lambda_2,\dots$ be variables. Then, we have
$$\sum_{\sigma\in S_{n+1}}\prod_{k=1}^n\frac{\lambda_{\sigma(1)}-t\lambda_{\sigma(n+1)}}{\lambda_{\sigma(k)}-\lambda_{\sigma(k+1)}}=\sum_{k=0}^n\binom{n}k^2t^k.$$
Let $A_n$ be the root lattice generated as a monoid by
$\{\pmb{e}_i - \pmb{e}_j: 0 \leq i,j\leq n+1, \text{and $i\neq j$}\}$. Let $P(A_n)$ be the polytope formed by the convex hull of this generating set. Then the coefficients of $\sum_k\binom{n}k^2t^k$ are the
$h$-vectors of a unimodular triangulation of $P(A_n)$. \rm
\smallskip
\noindent
\bf Proof. \rm Single out the values $\sigma(1)=i, \sigma(n+1)=j$ and apply Lemma 22 to obtain
$$\align \text{LHS}&=\sum\Sb i,j=1\\ i\neq j\endSb^{n+1}(\lambda_i-t\lambda_j)^n\sum_{S_{n-1}}
\frac1{(\lambda_i-\lambda_{\sigma(2)})(\lambda_{\sigma(n)}-\lambda_j)}\prod_{k=2}^{n-1}
\frac1{\lambda_{\sigma(k)}-\lambda_{\sigma(k+1)}} \\
&=-\sum\Sb i,j=1\\ i\neq j\endSb^{n+1}(\lambda_i-t\lambda_j)^n(\lambda_j-\lambda_i)^{n-2}
\prod\Sb k=1\\ k\neq i,j\endSb^{n+1}\frac1{(\lambda_k-\lambda_i)(\lambda_k-\lambda_j)}
=\sum_{i,j=1}^{n+1}
%(\lambda_i-t\lambda_j)^n(\lambda_j-\lambda_i)^n
\prod\Sb k=1\\ k\neq i\endSb^{n+1}\frac{\lambda_j-\lambda_i}{\lambda_k-\lambda_i}
\prod\Sb k=1\\ k\neq j\endSb^{n+1}\frac{\lambda_i-t\lambda_j}{\lambda_k-\lambda_j}. \endalign$$
We now invoke Proposition 0 to utilize invariance and replace $\lambda_k=k$. Some simplification leads to
$$\align \text{LHS}&=\frac1{n!^2}\sum_{i,j=0}^n(-1)^{i+j}\binom{n}i\binom{n}j[i+1-t(j+1)]^n(j-i)^n \\
%&=\frac1{n!^2}\sum_{i,j,k=0}^n(-1)^{i+j+k}\binom{n}i\binom{n}j\binom{n}kt^k(j+1)^k(i+1)^{n-k}(j-i)^n \\
&=\sum_{k=0}^n\binom{n}kt^k\left[\frac1{n!^2}
\sum_{i,j=0}^n(-1)^{i+j+k}\binom{n}i\binom{n}j(j+1)^k(i+1)^{n-k}(j-i)^n\right].\endalign$$
Denote the double-sum by $W(n,k)$. We prove $W(n-1,k)+W(n-1,k-1)=W(n,k)$, the Pascal recurrence:
$$\align W(n-1,k)+W(n-1,k-1)&=\sum_{i,j\geq0}(-1)^{i+j+k}\frac{(j+1)^{k-1}(i+1)^{n-1-k}(j-i)^{n-1}}{i!j!(n-i-1)!(n-j-1)!}
\left[(j+1)-(i+1)\right] \\
&=\sum_{i,j\geq0}(-1)^{i+j+k}\frac{(j+1)^{k}(i+1)^{n-k}(j-i)^n}{(i+1)!(j+1)!(n-i-1)!(n-j-1)!} \\
&=\sum_{i,j\geq1}(-1)^{i+j+k}\frac{j^{k}i^{n-k}(j-i)^n}{i!j!(n-i)!(n-j)!} \\
&=\sum_{i,j\geq0}(-1)^{i+j+k}\frac{j^{k}i^{n-k}(j-i)^n}{i!j!(n-i)!(n-j)!}=W(n,k);\endalign$$
where in the last equality once more we took the liberty (due to Proposition 0) of substituting $\lambda_k=k+z$, for a free parameter $z$.
In the present case, $z=-1$ is chosen. Together with the matching (check!) initial conditions, we conclude that $W(n,k)=\binom{n}k$ and hence the proof is complete. $\square$
\smallskip
\noindent
\bf Corollary 24. \it Let $\lambda_1,\lambda_2,\dots$ be variables. Taking $t=\pm1$ in Lemma 23 implies that
$$\sum_{\sigma\in S_{n+1}}\prod_{k=1}^n\frac{\lambda_{\sigma(1)}-\lambda_{\sigma(n+1)}}{\lambda_{\sigma(k)}-\lambda_{\sigma(k+1)}}=\binom{2n}n, \qquad
\sum_{\sigma\in S_{n+1}}\prod_{k=1}^n\frac{\lambda_{\sigma(1)}+\lambda_{\sigma(n+1)}}{\lambda_{\sigma(k)}-\lambda_{\sigma(k+1)}}=(-1)^{\binom{n+1}2}\delta_E(n)\binom{n}{n/2}.$$ \rm

\bigskip
\noindent
\bf PROBLEMS \rm

\bigskip
\noindent
\bf Problem 1. \it  Denote the Eulerian numbers by $A(n,j)$. Let $\lambda_1,\lambda_2,\dots$ be variables. Then, we have
$$\sum_{\sigma\in S_{n+1}}\prod_{k=1}^n\frac{\lambda_{\sigma(1)}+\cdots+\lambda_{\sigma(j)}}{\lambda_{\sigma(k)}-\lambda_{\sigma(k+1)}}=A(n,j-1)=\sum_{i=0}^j(-1)^i(j-i)^n\binom{n+1}i.$$ \rm

\bigskip
\noindent
\bf  Problem 2. \it  Let $Y, \lambda_1,\lambda_2,\dots$ be variables. Then, we have
$$\sum_{\sigma\in S_{n+1}}\prod_{k=1}^n\frac{\lambda_{\sigma(1)}^2+\cdots+\lambda_{\sigma(k)}^2-Y}{\lambda_{\sigma(k)}-\lambda_{\sigma(k+1)}}=(\lambda_1+\cdots+\lambda_{n+1})^n.$$ \rm

\bigskip
\noindent
\bf  Problem 3. \it Let $Y, \lambda_1,\lambda_2,\dots$ be variables, $p$ a positive integer and $s_{\mu}$ the Schur function for the partition $\mu=(p-1,\dots,p-1)=((p-1)^n)$. Then,
$$\sum_{\sigma\in S_{n+1}}\prod_{k=1}^n\frac{\lambda_{\sigma(k)}^p}{\lambda_{\sigma(k)}-\lambda_{\sigma(k+1)}}
=s_{\mu}(\lambda_1,\dots,\lambda_{n+1}).$$ \rm

\bigskip
\noindent
\bf  Problem 4. \it Let $Y, \lambda_1,\lambda_2,\dots$ be variables, $p, j$ positive integers and $m_{\mu}$ the monomial symmetric function of a partition $\mu$. Then, for $1\leq j\leq n+1$, we have
$$\sum_{\sigma\in S_{n+1}}\prod_{k=1}^n\frac{\lambda_{\sigma(j)}^p}{\lambda_{\sigma(k)}-\lambda_{\sigma(k+1)}}
=(-1)^{j-1}\binom{n}{j-1}\sum_{\mu\vdash (p-1)n}m_{\mu}(\lambda_1,\dots,\lambda_{n+1}).$$ \rm

\bigskip
\noindent
\bf  Problem 5. \it Let $Y, \lambda_1,\lambda_2,\dots$ be variables and $m_{\mu}$ is the monomial symmetric function associated with a partition $\mu$. Then, for $1\leq j\leq n+1$, we have
$$\sum_{\sigma\in S_{n+1}}\prod_{k=1}^n\frac{\lambda_{\sigma(1)}\cdots\lambda_{\sigma(j)}}{\lambda_{\sigma(k)}-\lambda_{\sigma(k+1)}}
=\sum\Sb \mu\vdash (j-1)n\\ \ell(\mu)\geq j\endSb m_{\mu}(\lambda_1,\dots,\lambda_{n+1}).$$ \rm

\bigskip
\noindent
\bf  Problem 6. \it Let $Y, \lambda_1,\lambda_2,\dots$ be variables, $p$ a positive integer and $m_{\mu}$ the monomial symmetric function of a partition $\mu$. Then, we have
$$\sum_{\sigma\in S_{n+1}}\prod_{k=1}^n\frac{\lambda_{\sigma(1)}^p\cdots\lambda_{\sigma(k)}^p}{\lambda_{\sigma(k)}-\lambda_{\sigma(k+1)}}
=\sum\Sb \mu\vdash p\binom{n+1}2-n\\ \mu\leq(pn-1,\dots,2p-1,p-1) \endSb
m_{\mu}(\lambda_1,\dots,\lambda_{n+1}).$$ \rm

\bigskip
\noindent
The following are certain generalizations to some of the formulas from the preceding pages, and as such it should not be difficult to emulate the techniques of this paper for the required justifications.
\bigskip
\noindent
\bf Exercise 7. \it Let $Y, \lambda_1,\lambda_2,\dots$ be variables. Then, we have
$$\align
&\sum_{\sigma\in S_{n+1}}\prod_{k=1}^n\frac{\lambda_{\sigma(k)}-Y_k}{\lambda_{\sigma(k)}-\lambda_{\sigma(k+1)}}
=1, \\
&\sum_{\sigma\in S_{n+1}}\prod_{k=1}^n\frac{1+\lambda_{\sigma(k)}Y_k}{\lambda_{\sigma(k)}-\lambda_{\sigma(k+1)}}
=Y_1Y_2\cdots Y_n, \\
&\sum_{\sigma\in S_{n+1}}\prod_{k=1}^n\frac{\lambda_{\sigma(1)}+\cdots+\lambda_{\sigma(k)}-Y_k}{\lambda_{\sigma(k)}-\lambda_{\sigma(k+1)}}=n!, \\
&\sum_{\sigma\in S_{n+1}}\prod_{k=1}^n\frac{(\lambda_{\sigma(1)}+\cdots+\lambda_{\sigma(k)})Y_k+1}{\lambda_{\sigma(k)}-\lambda_{\sigma(k+1)}}=n!Y_1Y_2\cdots Y_n, \\
&\sum_{\sigma\in S_{n+1}}
\frac{\sum_{k=1}^n(\lambda_{\sigma(1)}+\cdots+\lambda_{\sigma(k)})x_k-Y}
{\lambda_{\sigma(n)}-\lambda_{\sigma(n+1)}}
\prod_{k=1}^{n-1}\frac{\lambda_{\sigma(1)}+\cdots+\lambda_{\sigma(k)}-Y}{\lambda_{\sigma(k)}-\lambda_{\sigma(k+1)}}=(n-1)!\sum_{k=1}^nkx_k, \\
&\sum_{\sigma\in S_{n+1}}\prod_{k=1}^n\frac{\lambda_{\sigma(1)}-t\lambda_{\sigma(2)}}{\lambda_{\sigma(k)}-\lambda_{\sigma(k+1)}}=(1-t)^n-(n+1)(-t)^n, \\
&\sum_{\sigma\in S_{n+1}}\prod_{k=1}^n\frac{\lambda_{\sigma(1)}-t\lambda_{\sigma(3)}}{\lambda_{\sigma(k)}-\lambda_{\sigma(k+1)}}=(1-t)^n+[\binom{n}2-1](-t)^n-n\binom{n+1}2(-t)^{n-1}.
\endalign$$ \rm

\smallskip
\noindent
\bf QUESTION. \rm Is there a combinatorial or geometrical interpretation (as in [1]) to any of the formulas?


\Refs
\widestnumber\key{1}

\ref \key 1 \by Postnikov, A. \paper Permutohedra, associahedra, and beyond,
\jour Inter. Math. Res. Not. \vol 6 \#1-3 \yr2009 \pages1026-1106
\endref

\endRefs
\enddocument